\documentclass{ajmart}

\issueinfo{4}{4}{December}{2010}

\copyrightinfo{2007}{Aulona Press  \emph{(Albanian J. Math.)}}

\def\issn{{\sc ISSN} 1930-1235: }
\def\issueyear{2010}

\PII{\issn  (\issueyear )}

\pagespan{147}{160} 

\usepackage{amsmath,amsbsy,amsfonts,amsopn,amstext}
\usepackage{euscript,amssymb,amsbsy,amsfonts,amsthm,latexsym,amsopn,amstext, amsxtra,euscript,amscd}
\usepackage{graphicx, verbatim}

\usepackage{amsrefs}
\RequirePackage{rkeyval}[2001/12/22]
\RequirePackage{textcmds}[2001/12/14]
\RequirePackage{mathscinet}[2002/01/01]

\usepackage[all]{xy}
\input xy \xyoption{all}

\newtheorem{thm}{Theorem}

\newtheorem{lem}{Lemma}
\newtheorem{remark}{Remark}
\newtheorem{cor}{Corollary}
\theoremstyle{definition}

\theoremstyle{remark}

\def\x{x}
\def\y{y}
\def\z{z}

\def\bP{\mathbb P}

\def\X{\mathcal X}

\def\M{\mathcal M}

\def\T{\mathcal T}

\def\R{\mathbb R}

\def\P{{\mathbb P}}

\def\H_2{{\mathcal H_2}}  
\def\M{{\mathcal M}}

\def\S{\mathfrak S}

\def\e{\eta}
\def\G{\Gamma}

\def\lar{\longrightarrow}

\def\({\left(}
\def\){\right)}
\def\cO{{\mathcal O}}

\def\e{\varepsilon}

\def\P{\mathcal P}
\def\r{\delta}  
\def\<{\langle}
\def\>{\rangle}

\def\th{{\theta}}

\def\e{{\xi}}
\def\r{{r}}

\def\P{\mathcal P}
\def\<{\langle}
\def\>{\rangle}
\def\_u{{\mathfrak u}}
\def\J{J_{48}}

\def\det{\mbox{det }}

\def\J{\mbox{Jac }}

\def\H{\mathcal H}
\def\R{\mathbb R}

\def\P{{\mathbb P}}

\def\G2#1{{\Gamma^{(2)}#1}}
\def\Gs2#1{{\Gamma_0^{(2)}#1}}
\def\Ys2#1{{Y_0^{(2)}#1}}

\def\P{{\mathbb P}}

\begin{document}

\title[Singular locus  of Shaska surface]{Singular locus on the space of genus 2 curves with decomposable Jacobians.}

\author{Lubjana  Beshaj}

\address{Lubjana Beshaj, Department of Mathematics, University of Vlora,  Albania.}

\email{lbeshaj@univlora.edu.al}

\begin{abstract}
We study the singular locus on the algebraic surface $\S_n$ of genus 2 curves with a $(n, n)$-split Jacobian.  Such surface was computed by Shaska in \cite{deg3} for $n=3$,  and Shaska at al. in \cite{deg5} for $n=5$. We show that the singular locus for $n=2$ is exactly th locus of the curves of automorphism group $D_4$ or $D_6$. For $n=3$ we use a birational parametrization of the surface $\S_3$ discovered in \cite{deg3} to show that the singular locus is a 0-dimensional subvariety consisting exactly of three genus 2 curves (up to isomorphism) which have automorphism group $D_4$ or $D_6$.  We further show that the birational parametrization used in $\S_3$ would work for all $n \geq 7$ if $\S_n$ is a rational surface. 
\end{abstract}

\subjclass[2010]{14Q15, 14Q05, 68W30}

\keywords{genus two curves, moduli spaces, hyperelliptic curve cryptography, modular polynomials}

\setcounter{page}{147}

\maketitle


\begin{center}
\end{center}


\section{Introduction}
We study the singular locus on the space of genus 2 curves with a $(n, n)$-split Jacobian.  Such curves have been of much interest lately because of their use in many theoretical and applicative situations. The first part of the paper is based on several papers on the topic of genus two curves with split Jacobians; see 
\cite{ deg5, kyoto, neime, previato, sanjeewa, sh_u, sh_u2,   deg4, wijesiri, nato_wijesiri, caleb_wijesiri, sh_04, sh_03,   ajm_sh1, sh_02, sh_01, sh_05, nato_beshaj} among others.

In the first section, we study  genus 2 curves with split Jacobian. Let $\X$ be a genus 2 curve defined over an algebraically closed field $k$, of characteristic zero. Let $  \psi: \X \to  E $ be a degree $n$ maximal covering (i.e. does not factor through an isogeny) to an elliptic curve $E$ defined over $k$. We say that $\X$ has a degree $n$ elliptic subcover. Degree $n$ elliptic subcovers occur in pairs. Let $(E;E^{\prime})$  be such a pair. It is well known that there is an isogeny of degree $n^2$ between the Jacobian $\J(\X)$ of $\X$ and the product $E  \times E^{\prime}$. We say that $\X$ has $(n,n)$-split Jacobian.  

The locus of genus two curves with $(n, n)$-split Jacobians is an irreducible 2-dimensional algebraic variety.  There are many descriptions of it in the literature, but throughout this paper we will use only the embdedding of such space in the moduli space $\M_2$.  In other words, we would like an equation of such space where every point corresponds precisely to one isomorphism class of genus 2 curves. We denote such surface  by $\S_n$ and always think of it given by an equation in terms of the absolute invariants $i_1, i_2, i_3$ of genus two curves; see \cite{sh_05}. We will call the surface $\S_n$ the Shaska surface of level $n$.

The case with $(3, 3)$-split Jacobian was studied in \cite{deg3}.  These are the curves with degree $3$ elliptic subcovers. Shaska in ~\cite{deg3} computed the locus of curves $\X$ with degree 3 elliptic subfield in the moduli space of genus 2 curves.  We will give the explicit equation of this space and also  a graphical representation of it.  It was the first time that such an equation was computed other than the computationally trivial case for $n=2$. 

In \cite{deg5} was studied the case with $(5, 5)$ - split Jacobian by Shaska, Magaard, and Voelklein. There was  computed a normal form for the curves in the locus $\S_5$ and its three distinguished subloci. Further, they have  computed the equation of the elliptic subcover in all cases, gave a birational parametrization of the subloci of $\S_5$ as subvarieties of $\M_2$ and classify all curves in these loci which have extra automorphisms.

In section 2 of this paper  we  compute  the singular locus, $\T_2$, of the space $\S_2$, and the singular locus $\T_3$ of the space $\S_3$. The definition of the singular locus depends on the parametrization of the surface.  For the case of $n=2$ we prove that the singular locus of $\S_2$ is exactly the locus of genus 2 curves with automorphism group $D_4$ or $D_6$.  This computations were done using Maple 14.   

If the surface $\S_n$ is rational then we show how to obtain a birational parametrization for $\S_n$ using the invariants of binary cubics, which were used first in \cite{deg3}.

Throughout this paper by a genus two curve we mean the isomorphism class of a genus two curve defined over an algebraically closed field $k$. While most of the results are true for most characteristics, we assume throughout that the characteristic of $k$ is zero.

\section{Preliminaries}

\subsection{Genus 2 curves with split Jacobian}

Let $\X$ be a genus 2 curve defined over an algebraically closed field $k$, of characteristic zero. The affine version of this curve is given by the equation $\X: y^2 = F(x)$, where $F(x)$ is a polynomial of degree 5 or 6 and discriminant different from zero. Let 
\[\psi: \X  \rightarrow E\]
 be  a degree $n$ covering, where $n$ is odd and $E$ is an elliptic curve. The degree $n$ covering $\psi: \X  \rightarrow E$ induces a degree $n$ cover $\phi: \P^1 \to \P^1$  such that the following diagram commutes.
\[
\xymatrix {
  &   \X   \ar@{->}[dr]^{\pi_1} \ar@{->}[dl]_\psi  &          \\
%
%
E \ar@{->}[dr]_{\pi_2} &     &    \bP^1 \ar@{->}[dl]^\phi  \\
  &  \bP^1 & \\
}
\]
Here, $\pi_1 : \X \to \P_1$ and $\pi_2 : E \to \P_2$ are the hyperelliptic projections. So,  $\phi \circ \pi_1 = \pi_2 \circ \psi$. From Riemann- Hurwitz formula the number of  branch points is 4, or 5.   The ramification of the function $\phi$ is as follows; there are $\frac{n-1}{2}$ points of index 2 in $q_1$, $q_2$ and $q_3$, and  $\frac{n -3}{2}$ points of index 2 in $q_4$, and there is only one point of index 2 in $q_5$.  We denote this type of ramification by
$$ \left( (2)^{\frac{n-1}{2}}, (2)^{\frac{n-1}{2}}, (2)^{\frac{n-1}{2}}, (2)^{\frac{n-3}{2}}, (2) \right).$$
In the following figure bullets (resp., circles) represent  places of ramification index 2 (resp., 1).
\begin{figure}[!hbt]
$$
\xymatrix{& & & \circ \, \alpha_6  \\
\circ \, \alpha_1   & \circ  \, \alpha_2  & \circ \, \alpha_3  & \circ \, \alpha_5  & \circ \infty  \\
\bullet \, \ar[d] &\bullet \, \ar[d] &\bullet \, \ar[d] & \circ \, \alpha_4
  \ar[d]  &\bullet \,  0  \ar[d] \\
q_1 &q_2 &q_3 &  \infty  & 0      \\
}
$$
\caption{Ramification of $\phi : \P^1 \to \bP^1$ when $n=3$} \label{fig2}
\end{figure}
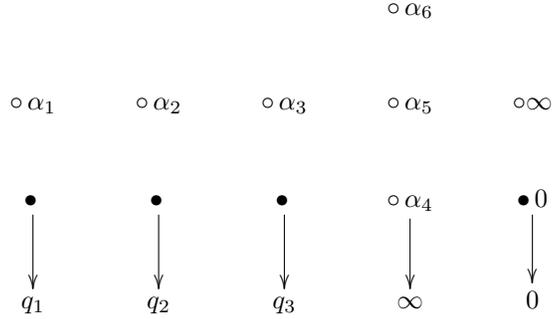

The family of coverings $\phi: \P^1 \to \P^1$, is an irreducible 2-dimensional algebraic variety. For every $\phi$ there exists a genus 2 curve $C$. Let $\H$ be the family of coverings.  We have the map
$$
\begin{aligned}
\alpha: \H & \to \M_2\\
[\phi] & \to [\X]
\end{aligned}
$$

Let $\alpha(\H)$ be denoted by $\S_n$. So, we say that these curves $\X$ are parametrized by an irreducible 2-dimensional subvariety $\S_n$ of
the moduli space $\M_2$ of genus 2 curves.  The fact that $\S_n$ is irreducible, for $n$ odd, comes from the braid action on Nielsen classes.  It is known that this is the case for all $n \cong 1 \mod 2$; see \cite{sh_01} among others.  Computation of spaces $\S_n$ as a subvariety of $\M_2$ has first computed by Shaska in \cite{deg3} for $n=3$ and then by Shaska, Magaard, and Voelklein for $n=5$; see \cite{deg5}.  We will call the  space $\alpha(\H)  \hookrightarrow \M_2$     the \textbf{Shaska surface of level $n$}.

\subsection{Pairs of elliptic subcovers}

Let $\psi_1:{\X} \lar E_1$ be a covering of degree $n$ from a curve of genus 2 to an elliptic curve.  The
covering $\psi_1:{\X} \lar E_1$ is called a \textbf{maximal covering} if it does not factor over a nontrivial
isogeny. A map of algebraic curves $f: X \to Y$ induces maps between their Jacobians $f^*: J_Y \to J_X$ and
$f_*: J_X \to J_Y$. When $f$ is maximal then $f^*$ is injective and $ker (f_*)$ is connected, see \cite{sh_01}
for details.

Let $\psi_1:{\X} \lar E_1$ be a covering as above which is maximal. Then ${\psi^*}_1: E_1 \to J_C$ is injective
and the kernel of $\psi_{1,*}: J_{\X} \to E_1$ is an elliptic curve which we denote by $E_2$, see \cite{sh_03} or
\cite{sh_05}. For a fixed Weierstrass point $P \in C$,  we can embed $C$ to its Jacobian via
$$i_P: {\X} \lar J_C$$
$$x \to [(x)-(P)]$$

Let $g: E_2 \to J_C$ be the natural embedding of $E_2$ in $J_C$, then there exists $g_*: J_{\X} \to E_2$.
Define $\psi_2=g_*\circ i_P: {\X} \to E_2$. So we have the following exact sequence
$$0 \to E_2 \buildrel{g}\over\lar J_{\X} \buildrel{\psi_{1,*}}\over\lar
E_1 \to 0$$ The dual sequence is also exact, see \cite{sh_01}
$$0 \to E_1 \buildrel{\psi_1^*}\over\lar J_{\X} \buildrel{g_*}\over\lar
E_2 \to 0$$
\par 

The following lemma shows that $\psi_2$ has the same degree as
$\psi_1$ and is maximal.
\begin{lem}
a) $deg\, (\psi_2) =n$

b) $\psi_2 $ is maximal
\end{lem}
For the proof see \cite{sh_01}.   If $deg (\psi_1)$ is an odd number then the maximal covering $\psi_2: {\X} \to E_2$ is unique (up to isomorphism of elliptic curves).

\par  To each of the covers $\psi_i:{\X} \lar E_i$, $i=1,2$, correspond 
covers $\phi_i: \bP^1 \lar \bP^1$. If the cover $\psi_1:{\X} \lar E_1$ is given, and therefore $\phi_1$, we want
to determine $\psi_2:{\X} \lar E_2$ and $\phi_2$. The study of the relation between the ramification structures
of $\phi_1$ and $\phi_2$ provides information in this direction.
The following lemma  answers this question for the set of Weierstrass points  $W=\{P_1, \dots , P_6\}$  of ${\X}$ when the degree of
the cover is odd.

\par Let $\psi_i:{\X} \lar E_i$, $i=1,2$, be maximal of odd degree $n$. Let ${\cO}_i\in E_i[2]$ be the points which has three Weierstrass points in its fiber.   Then, we have the following:

\begin{lem} The sets $\psi_1^{-1}({\cO }_1)\cap W$ and $\psi_2^{-1}({\cO }_2)\cap W$ form a disjoint union of W.
\end{lem}

Thus, the elliptic subcovers occur in pairs.
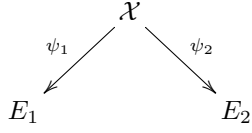
\begin{figure}
\[
\xymatrix {
  &   \X   \ar@{->}[dr]^{\psi_2}  \ar@{->}[dl]_{\psi_1}  &          \\
%
%
E_1   &     &    E_2    \\
%
%
}
\]
\caption{Splitting of the genus two curve}
\end{figure}

\subsection{Describing the Shaska surface  $\S_n$ in $\M_2$}

Consider a genus two curve $\X$ defined over $k$, given with equation
$$\X : \quad y^{2}=a_6 X^6+ a_5 X^5 +\dots  +a_0.$$
\emph{Igusa  $J$-invariants} $\, \, \{ J_{2i} \}$ of $\X$ are homogeneous polynomials of degree $2i$ in
\[k[a_0, \dots , a_6], \textit{  for }  i=1,2,3,5;\]
see \cite{sh_05}, \cite{sh_2000} for their definitions. Here $J_{10}$ is simply
the discriminant of $f(X,Z)$.   These
$J_{2i}$    are invariant under the natural action of $SL_2(k)$ on sextics. Dividing such an invariant by another
one of the same degree gives an invariant under $GL_2(k)$ action.

Two genus  2 fields $K$ (resp., curves) in the standard form $Y^2=f(X,1)$ are isomorphic if and only if the
corresponding sextics are $GL_2(k)$ conjugate. Thus if $I$ is a $GL_2(k)$ invariant (resp., homogeneous $SL_2(k)$
invariant), then the expression $I(K)$ (resp., the condition $I(K)=0$) is well defined. Thus the $GL_2(k)$
invariants are functions on the  moduli space $\mathcal M_2$ of genus 2 curves. This $\mathcal M_2$ is an affine
variety with coordinate ring
$$k[\mathcal M_2]=k[a_0, \dots , a_6, J_{10}^{-1}]^{GL_2(k)}$$
which is  the subring of degree 0 elements in $k[J_2, \dots ,J_{10}, J_{10}^{-1}]$. The \emph{ absolute
invariants}
$$ i_1:=144 \frac {J_4} {J_2^2}, \, \,  i_2:=- 1728 \frac {J_2J_4-3J_6} {J_2^3},\, \,  i_3 :=486 \frac {J_{10}}
{J_2^5}, $$
are even $GL_2(k)$-invariants. Two genus 2 curves with $J_2\neq 0$ are isomorphic if and only if they have
the same absolute invariants. If  $J_2=0 $ then we can define new invariants as in \cite{sh_05}. For the rest
of this paper if we say ``there is a genus 2 curve $\X$ defined over $k$'' we will mean the $k$-isomorphism
class of $\X$.

\begin{remark}
The definitions of $i_1, i_2, i_3$ with $J_{2}$ in the denominator is done simply for computational purposes. 
\end{remark}

Let 
\[F(X)=a_3X^3+a_2 X^2 + a_1 X+a_0, \textit{ and       }   G(X)=b_3 X^3+b_2X^2+b_1X+b_0\]
be two cubic polynomials.    We  define  the following invariants 
$$H(F,G) :=a_3 b_0 - \frac  1 3 a_2 b_1 + \frac 1 3 a_1 b_2 -a_0 b_3$$
We denote by $R(F,G)$ the resultant of $F$ and $G$ and by $D(F)$ the discriminant of $F$ always with respect to $X$. Also,
\[r_1(F,G) =\frac  {H(F,G)^3} {R(F,G)}, \quad r_2(F,G)=\frac {H(F,G)^4} {D(F)\, D(G)}. \]

In \cite{vishi} it is shown that $r_1, r_2$, and $ r_3=\frac {H(F,G)^2} {J_2(F\,G)}$ form a complete system
of invariants for unordered pairs of cubics.

Every curve $\X$ in $\S_n$ is written as a product of two cubics. In other words, its equation is 
\[ y^2= F(X) \cdot G(X)\]
for some $F(X), G(X) \in k[X]$.     We will use the invariants $r_1, r_2$ in relation with these cubics. Since the discriminants of such cubics can not be zero (otherwise the curve is not a genus two curve) then $D(F), D(G) $ are nonzero. For the same reason $F(X)$ and $G(X) $ don't have any common factors. Hence, $R(F, G) \neq 0$. Thus, $r_1, r_2$ are everywhere defined.

\section{Computation of singular locus $\T_{n}$}

Throughout this section we will use $x,y,z$ for absolute invariants $i_1, i_2, i_3$ respectively.  
 Let $\S_n$ be the Shaska  surface of level $n$ given by 
\[ \S_n (\x, \y, \z) =0\]

Then, its singular set is defined as the solution of the system 
\begin{equation}
\begin{split}
\left\{ \aligned
  \frac {\partial {\S_n}} {\partial  {\x}}=0\\
  \frac {\partial {\S_n}} {\partial  {\x}}=0\\
  \frac {\partial {\S_n}} {\partial  {\x}}=0\\
  \S_n(\x, \y, \z )=0\\
 \endaligned
\right.
\end{split}
\end{equation}

\subsection{The singular locus $\T_2$}  The equation of $\S_2$ is given by 
\begin{tiny}
\[
\begin{split}
\S_2 (x,y,z) &= -27\,x^6-9459597312000\,z^2\,x^2+20639121408000\,z^2\,y+111451255603200\,z^2\,x-240734712102912\, z^2\\
& -55240704\, z\, x^4-18\, y^2\, x^4  -8294400\, z\, y^2\, x^2-47278080\, z\, y\, x^3-264180754022400000\, z^3   \\
& -2866544640000\, z^2\, y\, x+2\, x^6\, y-4\, x^3\, y^3+9\, x^7 +331776\, z\, x^5+107495424\, z\, y\, x^2-27\, y^4+9\, x\, y^4\\
& -52254720\, z\, y^2\, x+2\, y^5+161243136\, z\, y^2+161243136\, z\, x^3 -12441600\, z\, y^3+54\, x^3\, y^2  =0 \\
\end{split}
\]
\end{tiny}

\begin{figure}[hbt]
\includegraphics[width=6.0cm]{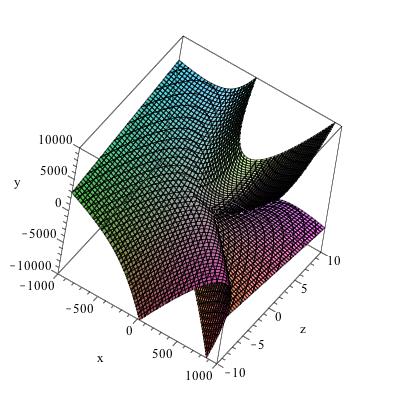}
\caption{The surface $\S_2$ graphed in $\R^3$.}
\end{figure}

Then we have the corresponding system from which we eliminate $z$ and get
\[ z=- {\frac {1}{82944}}\,  \frac {\phi_1(x, y)}  {\phi_2(x,y)}\]
where  $\phi_1$ and $\phi_2$ are as follows; 
\[
\begin{split}
\phi_1 (x,y) & = 104976\,{y}^{2}+5211\,{x}^{5}-48600\,{y}^{2}x+69984\,y{x}^{2}+3375\,y{x}^{4}+450\,{x}^{3}{y}^{2}\\
& -50544\,{x}^{4}-675\,{x}^{2}{y}^{2} +104976\,{x}^{3}+2025\,x{y}^{3}-10800\,{y}^{3}+20\,{x}^{6}+250\,{y}^{4}\\
& -37800\,{x}^{3}y\\
\phi_2 (x, y) & =  1250\,y{x}^{2}-121500\,xy - 3779136-359100\,{x}^{2} -11250\,{y}^{2} +6375\,{x}^{3} \\
& +421200\,y +2274480\,x \\ 
\end{split}
\]
The locus $\T_2$ which has 3 irreducible components which we describe below algebraically and graphically. 

\begin{figure}[hbt]
\includegraphics[width=6.0cm]{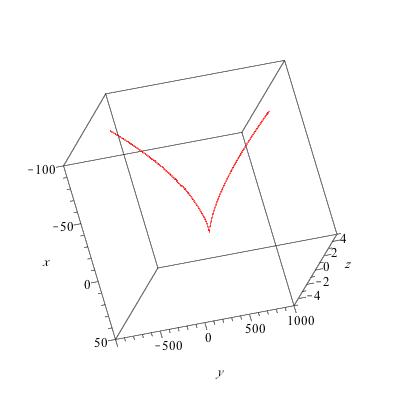}
\caption{The component $C_1$}
\end{figure}

The first component is given by 
\begin{small}
\[ C_1: \quad  100\,{y}^{2}-1458\,y+540\,xy-243\,{x}^{2}+80\,{x}^{3}    =0 \]
\end{small}
it corresponds to the locus of genus two curves with automorphism group $D_4$.

The second component is given by
\[  C_2:  3888\,x-1188\,{x}^{2}+5\,{x}^{3}+432\,y-360\,xy-25\,{y}^{2} =0 \]
and it corresponds to the locus of genus two curves with automorphism group $D_6$. 
\begin{figure}[hb]
\includegraphics[width=6.0cm]{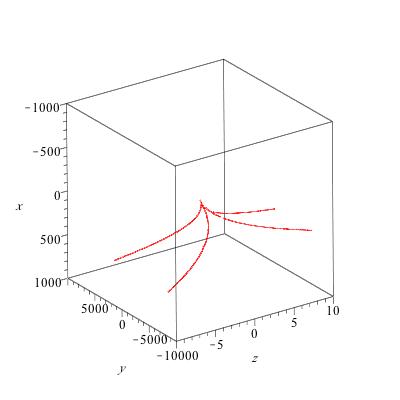}
\caption{The component $C_2$}
\end{figure}

The third component of $\T_2$ is given by the following system
\[
C_3: \quad \left\{
\begin{split}
& 50\,{x}^{4}-7515\,{x}^{3}-825\,y{x}^{2}+20412\,{x}^{2}-23490\,xy-4050  \,{y}^{2}+52488\,y =0\\
& 125\,{y}^{2}-1620\,y+1125\,xy-5832\,x+1890\,{x}^{2}+25\,{x}^{3} =0\\
\end{split}
\right.
\]

The solution of the $C_3$ system is 
\[
\left\{
\begin{split}
y & ={\frac {1}{75}}\,{\frac {408240\,x-33525\,{x}^{2}-944784+250\,{x}^{3}}{-864+55\,x}}\\
&  125\,{x}^{3}-9450\,{x}^{2}+247860\,x -944784   =0\\
\end{split}
\right.
\]
and the points $(x, y)$ given by 
\[ 
\begin{split}
\left(0, \frac {729} {50}\right), \left( \frac {81} {20}, - \frac {729} {200},   \right), 
\left( - \frac {36} 5, \frac {1512} {25}\right)
\end{split}
\]
However, only the first point is on the variety and it is \[ \left(0, \frac {729} {50}, \frac {729} {12800000} \right)\]
and has  automorphism groups are $D_4$ and therefore is contained in the first component.

We summarize  in the following theorem:

\begin{thm}
The singular locus of $\T_2$ contains two components, the irreducible loci of curves of automorphism group $D_4$ and $D_6$. 
\end{thm}

\subsection{The locus $\T_3$}

In this section we compute the singular locus $\T_3$ of $\S_3$.  The equation of $\S_3$ is quite large and was computed in \cite{deg3}. Below we display this equation $\S (x,y,z) \mod 5$. 

\medskip
\begin{tiny}
\[
\begin{split}
& x^{20}+3\,x^{19}+3\,x^{18}y+4\,x^{17}y^{2}+3\,x^{18}+4\,x^{17}z+2\,x^{16}y^{2}+2\,x^{16}yz+2\,x^{15}y^{3}+4\,x^{16}z  +2\,x^{15}y^{2}\\
& +4\,x^{15}yz+x^{15}{z}^{2}+x^{13}y^{3}z+3\,x^{14}yz+x^{13}y^{2}z+x^{13}y{z}^{2}+4\,x^{12}y^{3}z+4\,x^{12}y^{2}{z}^{2}+x^{11}y^{4}z+x^{10}y^{5}z\\
& +4\,x^{13}{z}^{2} + x^{12}y^{2}z+4\,x^{12}{z}^{3}+3\,x^{11}y^{3}z+3\,x^{11}y^{2}{z}^{2}+2\,x^{11}y{z}^{3}+4\,x^{10}y^{4}z+2\,x^{10}{y}^{3}{z}^{2}\\
& +2\,x^{9}y^{5}z+2\,x^{9}y^{4}{z}^{2}+2\,x^{8}{y}^{6}z+x^{7}y^{7}z+4\,x^{5}y^{10}+3\,x^{12}{z}^{2}+3\,x^{11}y{z}^{2}+3\,x^{11}{z}^{3}+4\,x^{10}y{z}^{3}+4\,x^{9}y^{4}z\\
& +3\,x^{9}y^{3}{z}^{2}+2\,x^{9}y^{2}{z}^{3}+3\,x^{8}y^{5}z+4\,x^{8}y^{4}{z}^{2}+3\,x^{8}y^{3}{z}^{3}+2\,x^{7}y^{6}z+2\,x^{7}y^{5}{z}^{2}+3\,x^{5}y^{8}z+2\,x^{4}y^{10}+x^{4}{y}^{9}z\\
& +2\,x^{3}y^{11}+x^{2}y^{12}+2\,x^{10}{z}^{3}+3\,x^{9}y^{2}{z}^{2}+4\,x^{9}y{z}^{3}+x^{9}{z}^{4}+4\,x^{8}y^{3}{z}^{2}+4\,x^{8}y^{2}{z}^{3}+2\,x^{8}y{z}^{4}+3\,x^{7}y^{4}{z}^{2}\\
& +2\,x^{6}y^{6}z+4\,x^{6}y^{5}{z}^{2}+2\,x^{6}y^{4}{z}^{3}+3\,x^{5}y^{7}z+x^{5}y^{5}{z}^{3}+4\,x^{4}y^{7}{z}^{2}+2\,x^{3}y^{10}+3\,x^{3}y^{9}z+4\,x^{3}y^{8}{z}^{2}+3\,xy^{12}\\
& +4\,xy^{11}z+3\,y^{13}+4\,x^{9}{z}^{3}+x^{8}y{z}^{3}+3\,{x}^{8}{z}^{4}+2\,x^{7}y^{2}{z}^{3}+2\,x^{7}y{z}^{4}+2\,x^{7}{z}^{5}+x^{6}y^{4}{z}^{2}+x^{6}y^{3}{z}^{3}+3\,x^{6}y^{2}{z}^{4}\\
& +x^{6}y{z}^{5}+4\,x^{5}y^{5}{z}^{2}+x^{5}y^{4}{z}^{3}+{x}^{5}y^{3}{z}^{4}+x^{4}y^{6}{z}^{2}+2\,x^{4}y^{5}{z}^{3}+x^{4}y^{4}{z}^{4}+3\,x^{3}y^{6}{z}^{3}+3\,x^{2}y^{9}z+3\,x^{2}y^{8}{z}^{2}\\
&+4\,x^{2}y^{7}{z}^{3}+4\,xy^{10}z+3\,y^{12}+2\,y^{11}z+x^{7}{z}^{4}+x^{6}y^{2}{z}^{3}+3\,x^{6}y{z}^{4}+3\,x^{6}{z}^{5}+4\,x^{5}y^{3}{z}^{3}+x^{5}y^{2}{z}^{4}+3\,x^{5}y{z}^{5}\\
&+3\,x^{5}{z}^{6}+2\,x^{4}y^{4}{z}^{3}+4\,x^{4}y^{3}{z}^{4}+x^{4}y^{2}{z}^{5}+4\,x^{3}y^{4}{z}^{4}+3\,x^{3}{y}^{3}{z}^{5}+2\,x^{2}y^{7}{z}^{2}+4\,x^{2}y^{6}{z}^{3}+2\,x^{2}y^{5}{z}^{4}\\
& +2\,xy^{8}{z}^{2}+3\,xy^{7}{z}^{3}+3\,y^{10}z+3\,y^{9}{z}^{2}+2\,x^{6}{z}^{4}+3\,x^{5}y{z}^{4}+3\,x^{5}{z}^{5
}+x^{4}y^{2}{z}^{4} +3\,x^{4}{z}^{6}+2\,x^{3}y^{3}{z}^{4}\\
&+3\,x^{3}y^{2}{z}^{5}+3\,x^{2}y^{5}{z}^{3}+3\,x^{2}y^{4}{z}^{4}+3\,xy^{6}{z}^{3}+2\,xy^{5}{z}^{4}+2\,xy^{4}{z}^{5}+2\,y^{7}{z}^{3}+y^{5}{z}^{5}+2\,x^{4}{z}^{5}+x^{3}y{z}^{5}\\
& +3\,x^{3}{z}^{6}+2\,x^{2}y^{3}{z}^{4}+2\,x^{2}y^{2}{z}^{5}+2\,x^{2}y{z}^{6}+2\,xy^{4}{z}^{4}+3\,y^{5}{z}^{4}+4\,y^{4}{z}^{5}+2\,x^{3}{z}^{5}+3\,x^{2}y{z}^{5}+4\,x^{2}{z}^{6}+xy^{2}{z}^{5}\\
&+3\,y^{2}{z}^{6}+x{z}^{6}+3\,y^{2}{z}^{5}+4\,{z}^{7}+3\,{z}^{6} =0
\end{split}
\]
\end{tiny}

\medskip

\noindent Let $\X$ be a genus 2 curve in the locus $\S_{3}$.  Then,  $\X$ is given by the equation
\begin{equation}
y^2 = (4x^3v^2+x^2v^2+2xv+1)(x^3v^2+x^2uv+xv+1),
\end{equation}
see \cite{sh_02} for details.   In \cite{deg3}   was computed the equation of $\S_3$ using the map
$$\theta: (u, v) \to (i_1, i_2, i_3)$$
where  the absolute invariants $i_1, i_2, i_3$ in terms of $u,v$ are
\begin{small}
\begin{equation}\label{eq_i}
\begin{split}
i_1 = &\frac {144}  {v(-405+252u+4u^2-54v-12uv+3v^2)^2}  (1188u^3-8424uv+u^4v-24u^4\\
& +14580v -66u^3v+138uv^2+297u^2v+945v^2-36v^3+9u^2v^2)  \\
i_2 =&- \frac {864} {v^2(-405+252u+4u^2-54v-12uv+3v^2)^3}
(-81v^3u^4+2u^6v^2+234u^5v^2\\
& +3162402uv^2-21384v^3u+26676v^4-473121v^3-72u^6v -5832v^4u+14850v^3u^2\\
& -72v^3u^3+324v^4u^2-650268u^3v-5940u^3v^2-3346110v^2+432u^6-1350u^4v^2\\
& +136080u^4v -7020u^5v-307638u^2v^2   \\
i_3 =& -243 \frac {(v-27)(4u^3-u^2v-18uv+4v^2+27v)^3}{v^3(-405+252u+4u^2-5
4v-12uv+3v^2)^5} \\
\end{split}
\end{equation}
\end{small}

The   map
$$\th: (u,v) \to (i_1, i_2, i_3)$$
given by \eqref{eq_i} which has degree 2 and it is defined when $J_2\neq 0$. For now we assume that $J_2\neq 0$
(The case $J_2=0$ is treated in Section 4.2, of \cite{deg3}). Denote the minors of the Jacobian matrix of $\th$ by
$M_1(u,v), M_2(u,v), M_3(u,v)$. The solutions of

\begin{equation}
\begin{split}
\left\{ \aligned
  M_1(u,v)= 0 \\
  M_2(u,v)=0\\
  M_3(u,v)=0\\
\endaligned
\right.
\end{split}
\end{equation}
consist of the (non-singular) curve
\begin{equation}
\begin{split}\label{n_3_iso1}
8v^3+27v^2-54uv^2-u^2v^2+108u^2v+4u^3v-108u^3=0\\
\end{split}
\end{equation}
and 7 isolated solutions which we display in Table 1,    together with the corresponding values
$(i_1, i_2, i_3)$, the automorphism group, and the number of elliptic subcovers.

\begin{small}
\begin{table}[!ht]
\renewcommand\arraystretch{1.5}
\noindent\[
\begin{array}{|c|c|c|c|c|}
\hline
(u,v) & (i_1, i_2, i_3) & Aut(K) & e_3(K) \\
\hline
(-\frac 7 2, 2) & J_{10}=0, \quad \text{no associated} &  &  \\
 & \text{genus 2 field K}  &  & \\
\hline
(-\frac {775} 8, \frac {125} {96}), & & & \\
 (\frac {25} 2, \frac {250} {9})&
- \frac {8019} {20}, -\frac {1240029} {200}, \frac {531441} {100000}     & D_4 & 2 \\
\hline
(27- \frac {77} 2 \sqrt{-1}, 23+ \frac {77} 9 \sqrt{-1}), &  & & \\
(27+ \frac {77} 2 \sqrt{-1}, 23- \frac {77} 9 \sqrt{-1}) & (\frac {729} {2116}, \frac {1240029} {97336},
\frac {531441}
{13181630464}  & D_4 & 2   \\
\hline
(-15+ \frac {35} 8 \sqrt{5}, \frac {25} 2 + \frac {35} 6 \sqrt{5}), &  & & \\
(-15- \frac {35} 8 \sqrt{5}, \frac {25} 2 - \frac {35} 6 \sqrt{5})&
81, - \frac {5103} {25}, -\frac {729} {12500}  &  D_6 & 2  \\
\hline
\end{array}
\]
\vspace{.3cm}
\caption{Exceptional points where $\det (Jac(\th))=0$} \label{tab1}
\end{table}
\end{small}

Notice that the curve given by Eq.~\eqref{n_3_iso1} corresponds to genus 2 curves with isomorphic degree 3 elliptic subcovers. Hence, the cover has  singular branch locus on such cases. We will see next how this can be avoided when we use the invariants of a pair of cubics. 


\subsection{Birational parametrization of $\S_3$}

For $F(X)=(4x^3v^2+x^2v^2+2xv+1)$ and $G(X)=(x^3v^2+x^2uv+xv+1)  $  
 we have
\begin{equation}
\begin{split}\label{eq_r}
r_1(F,G) &= 27\frac {v(v-9-2u)^3} {4v^2-18uv+27v-u^2v+4u^3}\\
r_2 (F,G)& = -1296 \frac  {v(v-9-2u)^4}   {(v-27)(4v^2-18uv+27v-u^2v+4u^3)}\\
\end{split}
\end{equation}

\begin{lem} The function field of $\S_3$ is given by    $k(\r_1, \r_2)$. In other words  $k(i_1, i_2, i_3)=k(\r_1, \r_2).$  
Moreover;
\begin{small}
\begin{equation}\label{i1_i2_i3}
\begin{split}
i_1 &= \frac 9 4 \frac {(13824\r_1^3\r_2^2+442368\r_1^2\r_2^3+5308416\r_1\r_2^4+192\r_1^4\r_2+\r_1^5+786432\r_1\r_2^3+9437184\r_2^4)}{\r_1(-1152\r_2^2+96\r_2\r_1+\r_1^2)^2} \\
i_2 &=\frac {27} {8\r_1^2(-1152\r_2^2+96\r_2\r_1+\r_1^2)^3 } (+79626240\r_1^4\r_2^4-4076863488\r_1^2\r_2^5+34560\r_1^6\r_2^2\\
& +12230590464\r_1^2\r_2^6+32614907904\r_1\r_2^6+14495514624\r_2^6 +288\r_1^7\r_2+2211840\r_1^5\r_2^3\\
& +\r_1^8-212336640\r_1^3\r_2^4+1528823808\r_1^3\r_2^5 -2359296\r_1^4\r_2^3)\\
i_3  &=-521838526464 \frac {\r_2^9} {\r_1^2(-1152\r_2^2+96\r_2\r_1+\r_1^2)^5}\\
\end{split}
\end{equation}
\end{small}
\end{lem}

\begin{figure}[hbt]
\includegraphics[width=8.0cm]{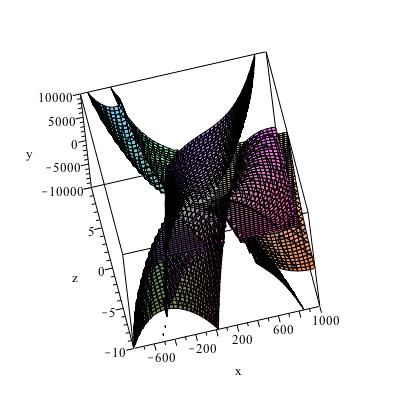}
\caption{Shaska surface  $\S_3$}
\end{figure}

The solution of the system in 
\begin{equation}
\begin{split}
\left\{ \aligned
  M_1(r_1, r_2)= 0 \\
  M_2(r_1, r_2)=0\\
  M_3(r_1, r_2)=0\\
\endaligned
\right.
\end{split}
\end{equation}
is
\begin{equation} \label{J2}
-1152 r_2^2+96 r_1 r_2+r_1^2 =0 
\end{equation}
and the system

\begin{tiny}
\[
\left\{
\begin{split}
& 3\,{r_1}^{8}+720\,{r_1}^{7}r_2+69120\,{r_1}^{6}{r_2}^{2}+2048\,{r_1}^{5}{r_2}^{2}+3317760\,{r_1}^{5}{r_2}^{3}+79626240\,{r_1}^{4}{r_2}^{4}-417792\,{r_1}^{4}{r_2}^{3} \\
& -24772608\,{r_1}^{3}{r_2}^{4} +764411904\,{r_1}^{3}{r_2}^{5}-113246208\,{r_1}^{2}{r_2}^{5}+50331648\,r_1{r_2}^{5}\\
& -5435817984\,r_1{r_2}^{6}-2415919104\,{r_2}^{6} =  0\\
& 9\,{r_1}^{5}+1296\,{r_1}^{4}r_2+62208\,{r_1}^{3}{r_2}^{2}-10240\,{r_1}^{2}{r_2}^{2}+995328\,{r_1}^{2}{r_2}^{3}+786432\,r_1{r_2}^{3}-2359296\,{r_2}^{4}  =0\\
& 9\,{r_1}^{8}+2160\,{r_1}^{7}r_2+207360\,{r_1}^{6}{r_2}^{2}+9953280\,{r_1}^{5}{r_2}^{3}+38912\,{r_1}^{5}{r_2}^{2}+238878720\,{r_1}^{4}{r_2}^{4}\\
& -3735552\,{r_1}^{4}{r_2}^{3}+2293235712\,{r_1}^{3}{r_2}^{5} -247726080\,{r_1}^{3}{r_2}^{4}+905969664\,{r_1}^{2}{r_2}^{5}\\
& +201326592\,r_1{r_2}^{5}-5435817984\,r_1{r_2}^{6}-4831838208\,{r_2}^{6}  =0 \\
\end{split}
\right.
\]
\end{tiny}

Then we get the following singular points
\[
(r_1, r_2)= \left(  -\frac {512} {2187},  - \frac {256}{6561} \right),     \left( \frac {2} {243},  \frac 1 {11664}\right), \left( - \frac {4000} {2187},  \frac {2500} {6561} \right)
\]
and the corresponding points (respectively) in $\S_{3}$ are:
\[
\begin{split}
(i_1, i_2, i_3)=
& \left( - \frac {8019}{20}, -\frac {1240029}{200},  - \frac {531441}{100000} \right) ,  \\
& \left(   81,  - \frac {5103}{25},  - \frac {729}{12500} \right),  \\
& \left( \frac {729}{2116},  \frac {1240029}{97336},  \frac {531441}{13181630464}   \right) \\
\end{split}
\]
which have automorphism groups respectively $D_4$, $D_4$, $D_6$, as seen from Table 1.

Notice that the Eq.~\eqref{J2} is exactly the case for $J_2=0$ where $i_1, i_2, i_3$ are not defined.  

\begin{cor}
The singular locus $\T_3$ of $\S_3$ are the points 
\begin{small}
\[
\begin{split}
 \left( - \frac {8019}{20}, -\frac {1240029}{200},  - \frac {531441}{100000} \right) ,  
 \left(   81,  - \frac {5103}{25},  - \frac {729}{12500} \right),   \left( \frac {729}{2116},  \frac {1240029}{97336},  \frac {531441}{13181630464}   \right) \\
\end{split}
\]
\end{small}
which have automorphisms group $D_4, D_4, D_6$ respectively.
\end{cor}

Notice that we have to use a parametrization in order to get the singular locus, because it is difficult computationally to compute this locus via partial derivatives.

\section{Some remarks for the general case.} Let's give a general approach how one can attempt to compute the surface $\S_n$ for $n \geq 7$. For $n\geq 7$ we get the first general case where the symmetries between the fourth and the fifth branch points which occur for degree 5 do not occur any longer; see \cite{deg5}. 

Suppose that $n \geq 7$.  Then $\S_n$ is parametrized by the $r_1, r_2$ invariants of two cubics. As in \cite{sh_01} we write a system of equations for the degree 7 covering $\phi: \P^1 \to \P^1$. 

Let $\X$ be a genus 2 curve in $\S_n$ which has equation 
\[ y^2 = (x^3+a x^2 + bx +c) \, (x^3+ux^2+vx +w)\]
such that  $a, b, c, u, v$ are expressed in terms of the two parameters $u$ and $v$. 
Let $r_1$ and $r_2$ be the invariants of the two cubics.  Then, 
there is a birational parametrization of $\S_n$ in terms of parameters $(r_1, r_2)$, i.e. 
\[ (r_1, r_2) \to (i_1, i_2, i_3 )\]
such that    $k(\S_n) = k(r_1, r_2)$.   Moreover, the singular locus of this parametrization contains the locus
\[ J_2 (r_1, r_2 )=0\]

While the computation of $\S_n$ for $n\geq 7$ is more difficult because the degree is larger, it is also true that there are no other symmetries now other than the $S_3$ action on the first three branch points as described in \cite{deg3} and \cite{deg5} for cases $n=3,5$ respectively.

\bigskip

\noindent \textbf{Acknowledgements:}   I would like to thank the Department of Mathematics at Oakland University for their support during the time that this article was written.

\end{document}